\documentclass[10pt]{amsart}

\usepackage{amscd}%
\usepackage{amsmath}%
\usepackage{amsfonts}%
\usepackage{amssymb}%
\usepackage{graphicx}
\usepackage[T2A]{fontenc}
\usepackage[russian,english]{babel}

\theoremstyle{plain}
\newtheorem{theorem}{Theorem}[section]
\newtheorem{lemma}{Lemma}[section]

\newtheorem{corollary}{Corollary}[section]
\newtheorem{remark}{Remark}[section]

\numberwithin{equation}{section}

\begin{document}

\title[On isometric representations of the perforated semigroup]{On isometric representations of the perforated semigroup}
\author{S.A.Grigoryan, V.H.Tepoyan}
\maketitle

\textbf{Abstract. }{  We study isometric representations of the semigroup $\mathbb{Z}_+\backslash \{1\}$. Notion of an 
 inverse representation is introduced and a complete description (up to unitary equivalence) of such representations is given. 
 
Also, we study a class of non-inverse irreducible representations -- $\beta$-representations of the semigroup
 $\mathbb{Z}_+\backslash \{1\}$.}

\textbf{KeyWords:}{ inverse semigroup, regular representation, isometric representation, inverse representation, $C^*$-algebras, $\beta$-representation. }

\section{Introduction}

 In \cite{Coburn} L.A.Coburn showed that all the isometric representation of the semigroup of non-negative integers generate canonically
 isomorphic $C^*$-algebras. Later in the papers \cite{Douglas} R.G.Douglas and G.J.Murphy \cite{Murphy} proved a similar statement for the 
semigroups with archimedean order and total order, respectively. In \cite{AT} it was shown: necessary and sufficient condition for all the  non-unitary isometric representations of the semigroup $S$ to generate canonically isomorphic $C^*$-algebras is that the natural order on $S$ is total. 

The simplest example of a semigroup with a non total order is a semigroup $\mathbb{Z}_+\backslash \{1\}$. 
For the first time this semigroup was mentioned in \cite{Murphy} by G.J. Murphy. S.Y. Jang in \cite{Jang2} pointed to the two representations
 of this semigroup that generate  canonically non-isomorphic $C^*$-algebras. I. Raeburn and S.T. Vittadello in \cite{Raeburn} have studied 
all the isometric representations of  $\mathbb{Z}_+\backslash \{1\}$ satisfying particular condition.

In this paper we study the isometric representations of the semigroup $\mathbb{Z}_+\backslash \{1\} $. We introduce the notion of an 
inverse representation and show that there exist only two (up to unitary equivalence) irreducible representations. It is proved that these
 are the representations which are contained in the papers \cite{Jang2}, \cite{Murphy}, \cite{Raeburn}, \cite{Vittadello}. Also we 
study the non-inverse representations of the semigroup $\mathbb{Z}_+\backslash \{1\}$.
\section{Inverse representations}\label{par2}

Let $S$ be an abelian additive cancellative semigroup with zero which does not contain a group other than the trivial one. 
We denote by $\Gamma$ the Grothendieck group generated by the semigroup $S$. Recall that the group $\Gamma$ is the quotient of the
 semigroup $S\times S$ by the equivalence relation: $(a,b)\sim(c,d)$ if and only if $a+d=b+c$. The inverse of the quotient class $[(a,b)]$
 is $[(b,a)]$. We write $\Gamma=S-S$.

Let $\pi:S\rightarrow B(H_{\pi})$ be an \emph{isometric} representation of the semigroup $S$ in the algebra $B(H_{\pi})$ of all bounded
 linear operators on a Hilbert space $H_{\pi}$. For each $a$ of $S$, in this paper, $ T_{\pi}(a)$ denotes the isometry $\pi(a)$.
$T^*_{\pi}(a)$ is an adjoint operator of $T_{\pi}(a)$. Thus, $T^*_{\pi}(a)T_{\pi}(a)=I$ is the identity operator and 
$T_{\pi}(a)T^*_{\pi}(a)=P_{\pi}(a)$ is projection ($P_{\pi}(a)\neq I$).

Operators $T_{\pi}(a)$ and $T^*_{\pi}(b)$, where $a,b\in S$ are called  \emph{trivial} monomials. We calle a \emph{monomial} the finite product of trivial monomials.

The set of all monomials form a multiplicative involutive semigroup which we denote by $S_{\pi}^*$.

Defineа order on $S$: $a\prec b$, if $b=a+c$. With respect to this order $S$ is a net.

\begin{lemma}
For any monomial $V$ there exist $a$ and $b$ in $S$ such that
$$\lim_{c\in{S}}{T^*_{\pi}(c)VT_{\pi}(c)}=T^*_{\pi}(a)T_{\pi}(b),$$
where $\lim_{c\in{S}}$ is a limit on the net $S$.
\end{lemma}

Note that if $T^*_{\pi}(a)T_{\pi}(b)=T^*_{\pi}(c)T_{\pi}(d)$ for some $a,b,c$ and $d$ of $S$, then $b+c=a+d$. Therefore, for each
 monomial $V$ we can associate a unique element $b-a$ of the Grothendieck group $\Gamma$. An element $b-a$ is called an \emph{index} of
 the monomial $V$ and denoted by $\textnormal{ind}V=b-a$. In fact, this proves the following statement.

\begin{lemma}
\ \ \ \ 
\begin{enumerate}
	\item $\textnormal{ind}V$ does not depend on the way of representing the monomial $V$ as a product of elementary monomials;
	\item $\textnormal{ind}V_1\cdot V_2=\textnormal{ind}V_1+\textnormal{ind}V_2$.
\end{enumerate}
\end{lemma}

We denote by $S_{0,\pi}^{*}$ a subsemigroup of the semigroup $S_{\pi}^{*}$ consisting of the monomials $V$ such that $\textnormal{ind}V=0$. We say
 that the isometric representation $\pi:S\rightarrow B(H_{\pi})$ is \emph{an inverse representation} if $S_{\pi}^{*}$ is an inverse
 semigroup under the operation of $\ast$-involution, or equialently, the semigroup $S_{0,\pi}^{*}$ is a semigroup of idempotents
 in $S_{\pi}^{*}$, that is, the semigroup of orthogonal projections. According to the Lemma~2.2 in \cite{Grigoryan} every semigroup $S$ has
 at least one inverse representation. On the other hand if the defined above order on $S$ is a total order, then all the isometric
 representations of the semigroup $S$ are inverse (see \cite{AT}).

 The simplest example of an inverse representation is a representation $\pi:\mathbb{Z}_+\rightarrow B(l^2(\mathbb{Z}_+))$ by shift 
operator $T_{\pi}(n)e_m=e_{n+m}$, where $e_n(m)=\delta_ {n,m}$ (Kronecker symbol), which form an orthonormal basis in
 $l^2(\mathbb{Z}_+)$. In this case the semigroup ${\mathbb{Z}_+}_{\pi}^{*}$ is a bicyclic semigroup.

\section{Inverse representations of the semigroup $\mathbb{Z}_+\backslash \{1\} $}\label{par3}

G.J.Murphy in \cite{Murphy}, and later S.Y.Jang in \cite{Jang2} have shown that for the semigroup $\mathbb{Z}_+\backslash \{1\}$ there exist at
 least two isometric representations which generate canonically non-isomorphic $C^*$-algebras. The isometric representations $\pi$ of the semigroup $\mathbb{Z}_+\backslash \{1\}$ when the projections $T_{\pi}(n)T^*_{\pi}(n)$ and 
$T_{\pi}(m)T^*_{\pi}(m)$ commute for any $n,m$ of $\mathbb{Z}_+\backslash\{1\}$, were investigated in \cite{Raeburn}. In this section we show that this condition implies 
that the representation $\pi$ is inverse. Here we obtain in fact the same results as in \cite{Raeburn}, but the proof of these results
 are based on the notion of \emph{ initial element } for the representations $\pi_i, \ i=1,2$. This paragraph is provided for full details.

Consider two isometric representations of the semigroup $\mathbb{Z}_+\backslash \{1\}$
$$\pi_1:\mathbb{Z}_+\backslash \{1\}\rightarrow B(l^2(\mathbb{Z}_+\backslash \{1\})), \ \ \pi_0:\mathbb{Z}_+\backslash \{1\}\rightarrow B(l^2(\mathbb{Z}_+))$$
defined by the shift operator
$$T_{\pi_1}(m)e_n=e_{m+n} \mbox{ and } \ T_{\pi_0}(m)f_n=f_{m+n},$$
where $\{e_n\}_{n\in\mathbb{Z}_+\backslash \{1\}}$ and $\{f_n\}_{n\in\mathbb{Z}_+}$ are natural orthonormal bases in
 $l^2(\mathbb{Z}_+\backslash \{1\})$ and $l^2(\mathbb{Z}_+)$ respectively. As  shown in \cite{Murphy} and \cite{Jang2} these representations are
 unitarily non equivalent. 

Representations $\pi_0$ and $\pi_1$ generate inverse semigroups $(\mathbb{Z}_+\backslash \{1\})_{\pi_0}^{*}$
 and $(\mathbb{Z}_+\backslash \{1\})_{\pi_1}^{*}$. Representation 
$\pi_1:\mathbb{Z}_+\backslash \{1\}\rightarrow B(l^2(\mathbb{Z}_+\backslash \{1\}))$ is inverse because it is regular (see \cite{AT}),
 and $\pi_0:\mathbb{Z}_+\backslash \{1\}\rightarrow B(l^2(\mathbb{Z}_+))$ is also inverse since 
$(T^*_{\pi_1}(2)T_{\pi_1}(3))^n=T_{\pi_1}(n)$ and $(\mathbb{Z}_+\backslash \{1\})_{\pi_1}^{*}$ is a bicyclic semigroup.

In this section we show that every pure inverse isometric (semi-unitary) irreducible representation of the semigroup
 $\mathbb{Z}_+\backslash \{1\}$ is unitarily equivalent either to $\pi_0$ or to $\pi_1$.

Let $\pi:\mathbb{Z}_+\backslash \{1\}\rightarrow B(H)$ be an isometric irreducible representation.

\begin{lemma}\label{lemma3.1}
$T^*_{\pi}(n)T_{\pi}(n+1)=T^*_{\pi}(n+1)T_{\pi}(n+2), \ n\neq 0.$
\end{lemma}

\begin{corollary}
$T^*_{\pi}(n)T_{\pi}(m)=T^*_{\pi}(n+l)T_{\pi}(m+l)$ for all $l$ of $\mathbb{Z}_+$.
\end{corollary}

For the further we need following obvious relations
$$T^*_{\pi}(2)=T^*_{\pi}(3)T^*_{\pi}(2)T_{\pi}(3), \ T^*_{\pi}(3)=T^*_{\pi}(4)T^*_{\pi}(2)T_{\pi}(3), \ T^*_{\pi}(3)T_{\pi}(2)=T^*_{\pi}(2)T^*_{\pi}(2)T_{\pi}(3).$$
From these relations we have immediately: if $h_0\in\ker{T^*_{\pi}(2)T_{\pi}(3)}$, then
$$T^*_{\pi}(2)T_{\pi}(3)h_0=T^*_{\pi}(2)h_0=T^*_{\pi}(3)h_0=T^*_{\pi}(3)T_{\pi}(2)h_0=0.$$

\begin{lemma}\label{lemma3.2}
Let $\pi:\mathbb{Z}_+\backslash \{1\}\rightarrow B(H)$ be an inverse irreducible representation. Suppose
 $\ker{T^*_{\pi}(2)T_{\pi}(3)}\neq{0}$. Then this representation is unitary equivalent to the representation
$\pi_1:\mathbb{Z}_+\backslash \{1\}\rightarrow B(l^2(\mathbb{Z}_+\backslash \{1\}))$.
\end{lemma}

Note that Lemma~\ref{lemma3.2} holds for any not necessarily inverse isometric representation of the semigroup 
$\mathbb{Z}_+\backslash \{1\}$.

Suppose now $\pi:\mathbb{Z}_+\backslash \{1\}\rightarrow B(H)$ is an inverse representation, that is, 
$(\mathbb{Z}_+\backslash \{1\})^*_{\pi}$ is an inverse semigroup. Then
$$P=T^*_{\pi}(3)T_{\pi}(2)T^*_{\pi}(2)T_{\pi}(3) \ \mbox{ and } \ Q=T^*_{\pi}(2)T_{\pi}(3)T^*_{\pi}(3)T_{\pi}(2)$$
are projections and $Q<P$.

Indeed
$$PQ=T^*_{\pi}(3)T_{\pi}(2)T^*_{\pi}(2)T_{\pi}(3)T^*_{\pi}(2)T_{\pi}(3)T^*_{\pi}(3)T_{\pi}(2)=$$
$$T^*_{\pi}(3)T_{\pi}(2)T^*_{\pi}(2)T_{\pi}(3)T^*_{\pi}(3)T_{\pi}(4)T^*_{\pi}(3)T_{\pi}(2)=$$
$$T^*_{\pi}(3)T_{\pi}(3)T^*_{\pi}(3)T_{\pi}(2)T^*_{\pi}(2)T_{\pi}(4)T^*_{\pi}(3)T_{\pi}(2)=$$
$$T^*_{\pi}(3)T_{\pi}(2)T^*_{\pi}(2)T_{\pi}(2)T_{\pi}(2)T^*_{\pi}(3)T_{\pi}(2)=$$
$$T^*_{\pi}(3)T_{\pi}(2)T_{\pi}(2)T^*_{\pi}(3)T_{\pi}(2)=T^*_{\pi}(3)T_{\pi}(4)T^*_{\pi}(3)T_{\pi}(2)=$$
$$T^*_{\pi}(2)T_{\pi}(3)T^*_{\pi}(3)T_{\pi}(2)=Q.$$

During the proof of this inequality we have used the Lemma~\ref{lemma3.1} and relations
$$T_{\pi}(2)T^*_{\pi}(2)T_{\pi}(3)T^*_{\pi}(3)=T_{\pi}(3)T^*_{\pi}(3)T_{\pi}(2)T^*_{\pi}(2), \ \ T^*_{\pi}(2)T_{\pi}(4)=T_{\pi}(2).$$

\begin{lemma}\label{lemma3.3}
Let $\pi:\mathbb{Z}_+\backslash \{1\}\rightarrow B(H)$ be an irreducible non unitary inverse representation of the semigroup
$\mathbb{Z}_+\backslash \{1\}$. Suppose that $\ker{T^*_{\pi}(2)T_{\pi}(3)}=\{0\}$. Then this representation is unitary equivalent to
 the representation
$$\pi_0:\mathbb{Z}_+\backslash \{1\}\rightarrow B(l^2(\mathbb{Z}_+)).$$
\end{lemma}

The next statement follows immediately from the Lemmas~\ref{lemma3.2} and \ref{lemma3.3}.

\begin{theorem}
Let $\pi:\mathbb{Z}_+\backslash \{1\}\rightarrow B(H)$ be a non-unitary isometric inverse representation. Then $\pi$ is unitarily
 equivalent either to $\pi_0$ or to $\pi_1$.
\end{theorem}

Thus we have the following theorem

\begin{theorem}
Every isometric inverse representation $\pi:\mathbb{Z}_+\backslash \{1\}\rightarrow B(H)$ can be represented as a direct sum
$$\pi=k\pi_0\oplus l\pi_1\oplus\pi_2,$$
where $k$ is a multiplicity of the representation $\pi_1:\mathbb{Z}_+\backslash \{1\}\rightarrow B(l^2(\mathbb{Z}_+\backslash \{1\}))$,
 $l$ is a multiplicity of the representation $\pi_0:\mathbb{Z}_+\backslash \{1\}\rightarrow B(l^2(\mathbb{Z}_+))$ and  $\pi_2$ is a 
unitary representation of the semigroup $\mathbb{Z}_+\backslash \{1\}$.
\end{theorem}

\begin{remark}
\end{remark}
It can be shown that for an isometric representation $\pi:\mathbb{Z}_+\backslash \{1\}\rightarrow B(H_{\pi})$ the inversion of the semigroup 
$(\mathbb{Z}_+\backslash \{1\})^*_{\pi}$ is equivalent to the commutativity of projections 
$P_{\pi}(n)=T_{\pi}(n)T^*_{\pi}(n), \ n\in\mathbb{Z}_+\backslash \{1\}$. The problem of equivalence of   the inversion and the commutativity of elementary projections for any isometric representation remains unsolved open.

\section{$\beta$-representation of the semigroup $\mathbb{Z}_+\backslash \{1\}$}\label{par4}

It follows immediately from Theorem of Coburn \cite{Coburn} that the semigroup $\mathbb{Z}_+$ admits only one, up to unitary equivalence,
 infinite irreducible representation. In this section we show that there exist cotinuum of such representations for a ''deformed'' semigroup 
$\mathbb{Z}_+\backslash \{1\}$.

Let $H_0$ be a Hilbert subspace of $l^2(\mathbb{Z}_+)$ generated by the basis $\{e_n\}_{n=2}^{\infty}, \ e_n(m)=\delta_{n,m}$. Denote
 by $H_{\beta}$ the Hilbert subspace of $l^2(\mathbb{Z}_+)$ generated by the elements of $H_0$ and the function $g_{\beta}=\beta e_0+te_1$,
 where $\beta\in\mathbb{C}, t\in\mathbb{R}_+$ and $\beta^2+t^2=1$. Obviously $g_{\beta}$ with the family $\{e_n\}_{n=2}^{\infty}$
is an orthonormal basis in $H_{\beta}$ ш $H_{\beta}=\mathbb{C}g_{\beta}\oplus H_0$.

Let $P_{\beta}:l^2(\mathbb{Z}_+)\rightarrow H_{\beta}$ be an orthogonal projection from $l^2(\mathbb{Z}_+)$ to $H_{\beta}$. Define the 
representation $\tau_{\beta}:\mathbb{Z}_+\backslash \{1\}\rightarrow B(H_{\beta})$ assuming that
$$\tau_{\beta}(n)=P_{\beta}\pi_{0}(n)P_{\beta},$$
where $\pi_0:\mathbb{Z}_+\backslash \{1\}\rightarrow B(l^2(\mathbb{Z}_+))$ is an inverse representation defined as above. Since
$T_{\pi_1}(n)$ is a mapping of $H_{\beta}$ to itself, the representation $\pi_{\beta}:\mathbb{Z}_+\backslash \{1\}\rightarrow B(H_{\beta})$
is an isometric representation of the semigroup $\mathbb{Z}_+\backslash \{1\}$.

\begin{lemma}
The representation $\tau_{\beta}:\mathbb{Z}_+\backslash \{1\}\rightarrow B(H_{\beta})$ is an inverse representation if and only if
 $\beta=0$ or $|\beta|=1$.
\end{lemma}

It  follows immediately from this Lemma, that all the representations $\tau_{e^{i\theta}}, \ 0<\theta\leq 2\pi$ are unitarily equivalent 
to the inverse representation $\pi_1:\mathbb{Z}_+\backslash \{1\}\rightarrow B(l^2(\mathbb{Z}_+\backslash \{1\}))$, and $\tau_0$ to 
the inverse representation $\pi_0:\mathbb{Z}_+\backslash \{1\}\rightarrow B(l^2(\mathbb{Z}_+))$.

\begin{theorem}
\ \
\begin{enumerate}
	\item\label{p1} $\tau_{\beta}:\mathbb{Z}_+\backslash \{1\}\rightarrow B(H_{\beta})$ is an irreducible isometric representation;
	\item\label{p2} If $\beta_1\neq\beta_2$, $|\beta_1|<1$, then the representations $\tau_{\beta_1}$ and $\tau_{\beta_2}$ are unitarily non equivalent.
\end{enumerate}
\end{theorem}

{\small \vspace{\baselineskip}\hrule \vspace{3pt}
\par
{\bf Grigoryan S.A.} -- {Kazan State Power Engineering University
\par
E-mail: {\it gsuren@inbox.ru} }
{\small \vspace{\baselineskip}\hrule \vspace{3pt}
\par
{\bf Tepoyan V.H.} -- Kazan State Power Engineering University
\par
E-mail: {\it tepoyan.math@gmail.com} }

\end{document}